\documentclass{commat}

\hypersetup{
	colorlinks=true,
	linkcolor=blue,
	filecolor=magenta,      
	urlcolor=cyan,
}
\newtheorem{question}{Question}
\newtheorem{conj}{Conjecture}

\newcommand{\C}{\mathbb{C}}
\newcommand{\R}{\mathbb{R}}
\newcommand{\Q}{\mathbb{Q}}
\newcommand{\K}{\mathcal{K}}
\newcommand{\Z}{\mathbb{Z}}
\newcommand{\PP}{\mathcal{P}}
\newcommand{\QQ}{\mathcal{Q}}

\newcommand{\Lat}{\mathcal{L}}

\newcommand{\conv}{\mathrm{conv}}

\title{%
	The integer point transform as a complete invariant
}

\author{%
	Sinai Robins
}

\authorinfo[%
Sinai Robins]{%
	Instituto de Matem\'atica e Estat\'\i stica, Universidade de S\~ao Paulo\\ 
	Rua do Mat\~ao 1010, 05508-090 S\~ao Paulo/SP, Brazil}{%
	srobins@ime.usp.br \\
	The author gratefully acknowledges the support of FAPESP 2020/08343-8
}

\abstract{%
	The integer point transform $\sigma_\PP$ is an important invariant of a rational polytope $\PP$, and here we show that it is a complete invariant.  We prove that it is only necessary to evaluate $\sigma_\PP$ at one algebraic point in order to uniquely determine $\PP$, by employing the Lindemann-Weierstrass theorem.  Similarly, we prove that it is only necessary to evaluate the Fourier transform of a rational polytope $\PP$ at a single algebraic point, in order to uniquely determine $\PP$. We prove that identical uniqueness results also hold for integer cones. 
	
	In addition, by relating the integer point transform to finite Fourier transforms, we  show that a finite number of \emph{integer point evaluations} of $\sigma_\PP$ suffice in order to uniquely determine $\PP$.  We also give an equivalent condition for central symmetry of a finite point set, in terms of the integer point transform, and prove some facts about its local maxima.  Most of the results are proven for arbitrary finite sets of integer points in $\R^d$.
}

\keywords{%
	Integer point transform, integer points, rational polytope, lattices, complete invariant, finite Fourier analysis, Lindemann-Weierstrass theorem
}

\msc{%
	AMS classification 52C07, 11H06, 11P21 
}

\VOLUME{31}
\NUMBER{2}
\YEAR{2023}
\firstpage{157}
\DOI{https://doi.org/10.46298/cm.11218}

\begin{document}

	\section{Introduction}\label{sec:intro}
	
	A polytope $\PP\subset \R^d$ is called an \emph{integer polytope}
	(respectively, a~\emph{rational polytope})
	if the coordinates of all its vertices are integers (respectively, rationals).
	Given a~rational polytope $\PP \subset \R^d$, we define its integer point transform
	by
	\begin{equation}\label{def: integer point transform}
		\sigma_{\PP}(\xi):= \sum_{n\in \Z^d \cap \PP} e^{2\pi i \langle n, \xi \rangle},
	\end{equation}
	for all $\xi \in \R^d$. Throughout, the word polytope refers to a~convex set, by definition.
	We observe that more generally, Definition~\eqref{def: integer point transform}
	makes sense for an arbitrary finite set of
	integer points $S\subset \Z^d$:
	\begin{equation}\label{more general integer point transform}
		\sigma_S(\xi):= \sum_{n\in S} e^{2\pi i \langle n, \xi \rangle}.
	\end{equation}
	Indeed, we will often develop general principles for any finite set of integer points in $\R^d$, and then later pass to polytopes, using the assumption of convexity.

	One of the main utilities of the integer point transform is the special evaluation at the origin:
	\begin{equation}\label{0 gives the number of integer points}
		\sigma_{\PP}(0):= \sum_{n\in \Z^d \cap \PP} 1 = \left | \Z^d \cap \PP \right |,
	\end{equation}
	the number of integer points in $\PP$. So we see that if $\PP$ is a~polytope, then the integer point transform discretizes the volume of
	$\PP$ in this sense. In fact, the integer point transform was so named because, by using a~lattice, it discretizes the Fourier transform of $\PP$, which is defined by
	$\int_\PP e^{2\pi i \langle \xi, x \rangle} dx$ (see~\cite{Robins}, for example).
	We already know from~\eqref{0 gives the number of integer points} that if we have any two rational polytopes $\PP, \QQ \subset \R^d$, then
	$\sigma_{\PP}(0) =\sigma_{\QQ}(0) \implies |\PP\cap \Z^d | = |\QQ\cap \Z^d |$. But is it possible that
	knowledge of $\sigma_{\PP}(\xi)$, for some finite collection of points $\xi\in \R^d$, might help us to uniquely identify $\PP$?
	This is our main motivating question.
	
	Historically, the importance of the integer point transform $\sigma_\PP$ for a~rational polytope $\PP$ surfaced naturally in combinatorial geometry, namely in Ehrhart's theory of integer point enumeration in polytopes (\cite{BeckRobins},~\cite{Robins}). In particular, Ehrhart's main theorem follows quickly from Brion's theorem, which enables us to write
	$\sigma_\PP$ as a~finite linear combination of exponential-rational functions~\cite{Robins} (see Remark~\ref{fourth remark}). The work of Fink, M\'esz\'aros, and Dizier~\cite{Fink.etal} gives applications of integer point transforms to Schubert polynomials. The recent work of
	Katharina Jochemko~\cite{Jochemko} shows that the sequence of integer point transforms
	$ {\sigma_{kP+Q}}\big|_{k\geq 0}$ satisfies a~multivariate linear recursion, where $\PP$ is an integer polytope, and $Q$ is any polytope.
	Here we do not assume knowledge of Ehrhart theory, and rather proceed from first principles.
	
	The integer point transform~\eqref{def: integer point transform}
	clearly lives on the torus - in other words, $\sigma_{\PP}$ is periodic on $\R^d$, with a~fundamental domain $[0, 1)^d$. We will therefore often restrict attention to the torus, or equivalently to
	the half-open cube $[0, 1)^d$.

	It is elementary that given any positive integer $N$,
	there are infinitely many distinct $d$-dimensional
	integer polytopes $\PP$ with $\left | \Z^d \cap \PP \right | = N$.
	Even after we mod out by the action of the
	modular group $GL_d(\Z)$, there are in general many (though finitely many) distinct integer polytopes $\PP$ with
	$\left | \Z^d \cap \PP \right | = N$.
	It is very natural to ask: ``What extra information do we need in order to
	uniquely determine $\PP$?"
	To make this question more rigorous, we formulate it as follows.
	\begin{question}\label{first question}
		Given any two integer polytopes $\PP, \QQ \subset \R^d$,
		is there a~\emph{finite set} $S\subset \R^d$ such that
		\[
		\sigma_{\PP}(\xi) =\sigma_{\QQ}(\xi) \text{ for all } \xi \in S \iff \PP = \QQ?
		\]
	\end{question}
	
	\begin{question}\label{second question}
		More generally, given any two finite sets of integer points $A, B \subset \Z^d$,
		is there a~\emph{finite set} $S\subset \R^d$ such that
		\[
		\sigma_{A}(\xi) =\sigma_{B}(\xi) \text{ for all } \xi \in S \iff A=B?
		\]
	\end{question}
	We will answer both of these Questions in the affirmative.
	Somewhat surprisingly, it turns out that just one point suffices for both questions.
	We also prove a~similar, but slightly weaker result, for rational polytopes (Theorem~\ref{first main theorem}, part~\ref{part 3 of main thm} below).
	Because the direction $(\Longleftarrow)$ is trivial, it suffices to always prove the $(\Longrightarrow)$ direction.
	Throughout the paper, we will use the special point
	\begin{equation}\label{alg. indep. vector}
		\xi^* := \tfrac{1}{\pi}\left(\sqrt 2, \dots, \sqrt p_d\right)\in \R^d,
	\end{equation}
	where we have picked the first $d$ primes $2, 3, 5, 7, \dots p_d$ to ensure that
	$\sqrt 2, \dots, \sqrt p_d$ are algebraically independent over $\Q$, and therefore all integer linear combinations of the coordinates of $\xi^*$ are distinct.

	\begin{theorem}\label{first main theorem}
		\begin{enumerate}
			\item
			\label{part 1 of main thm}
			Fix any two finite sets of integer points $A, B \subset \Z^d$. Then:
			\begin{equation*}
				\sigma_{A}(\xi^*) =\sigma_{B}(\xi^*) \iff A = B.
			\end{equation*}
			\item
			\label{part 2 of main thm}
			\label{second part of main theorem}
			Let $\PP, \QQ \subset \R^d$ be integer polytopes. Then:
			\begin{equation}\label{first implication}
				\sigma_{\PP}(\xi^*) =\sigma_{\QQ}(\xi^*) \iff \PP = \QQ.
			\end{equation}
			\item
			\label{part 3 of main thm}
			Let $\PP, \QQ \subset \R^d$ be rational polytopes. Then:
			\begin{equation}\label{second implication}
				\sigma_{k\PP}(\xi^*) =\sigma_{k\QQ}(\xi^*) \iff \PP = \QQ,
			\end{equation}
			for any $k\in \Z_{>0}$ such that $k\PP$ and $k\QQ$ are both
			integer polytopes.
		\end{enumerate}
	\end{theorem}
	\begin{proof}
		We suppose that $\sigma_{A}(\xi^*) =\sigma_{B}(\xi^*)$.
		Then we have:
		\begin{align}
			0 &=\sum_{n\in A} e^{2\pi i \langle n, \xi^* \rangle} -
			\sum_{n\in B} e^{2\pi i \langle n, \xi^* \rangle} \\
			&=\sum_{n\in A} e^{2i (n_1\sqrt 2 +\cdots + n_d\sqrt{p_d})}
			- \sum_{n\in B} e^{2i (n_1\sqrt 2 +\cdots + n_d\sqrt {p_d})}. \label{vanishing sum}
		\end{align}
		To prove part~\ref{part 1 of main thm}, suppose to the contrary, that $A\not=B$.
		Then~\eqref{vanishing sum}
		gives us
		a finite \emph{nontrivial} vanishing sum of an integer linear combination of exponentials, all of which have
		the form $e^{\alpha}$, with $\alpha$ algebraic. But this contradicts the celebrated
		Lindemann–Weierstrass theorem: if $\alpha_1, \dots, \alpha_n$
		are distinct algebraic numbers, then $e^{\alpha_1}, \dots, e^{\alpha_n}$ are algebraically independent over $\Q$
		(\cite{Baker}, Theorem 1.4).

		To prove part~\ref{part 2 of main thm}, suppose that $\PP$ and $\QQ$ are integer polytopes.
		Part~\ref{part 1 of main thm} applies to the two sets of integer points
		$A:= \QQ \cap \Z^d$ and $B:= \PP \cap \Z^d$,
		and therefore $\QQ \cap \Z^d = \PP \cap \Z^d$.
		By invoking the very simple Lemma~\ref{simple lemma, going from point set to polytope} below, we conclude that $\PP =\QQ$.
		
		To prove part~\ref{part 3 of main thm},
		suppose that $\PP$ and $\QQ$ are any rational polytopes that enjoy
		the hypothesis $\sigma_{k\PP}(\xi^*) =\sigma_{k\QQ}(\xi^*)$. The positive integer dilate $k$ allows us to conclude
		that $k\PP$ and $k\QQ$ are both integer polytopes. By part~\ref{part 2 of main thm},
		$\sigma_{k\PP}(\xi^*) =\sigma_{k\QQ}(\xi^*) \implies k\PP = k\QQ$, which in turn implies that
		$ \PP = \QQ$.
	\end{proof}
	We note that no assumption had to be made about the dimensions of $\PP$ or $\QQ$. Furthermore,
	it is important to remark that the proof shows that more is true:
	we may pick any algebraic numbers in the exponents of our exponentials in equation~\eqref{vanishing sum},
	so in choosing $\xi^*$, there is clearly a~dense set of vectors to choose from.

	
	\section{Passing from finite point sets to polytopes}
	
\begin{figure}[htb]
		\begin{center}
			\includegraphics[totalheight=2.9in]{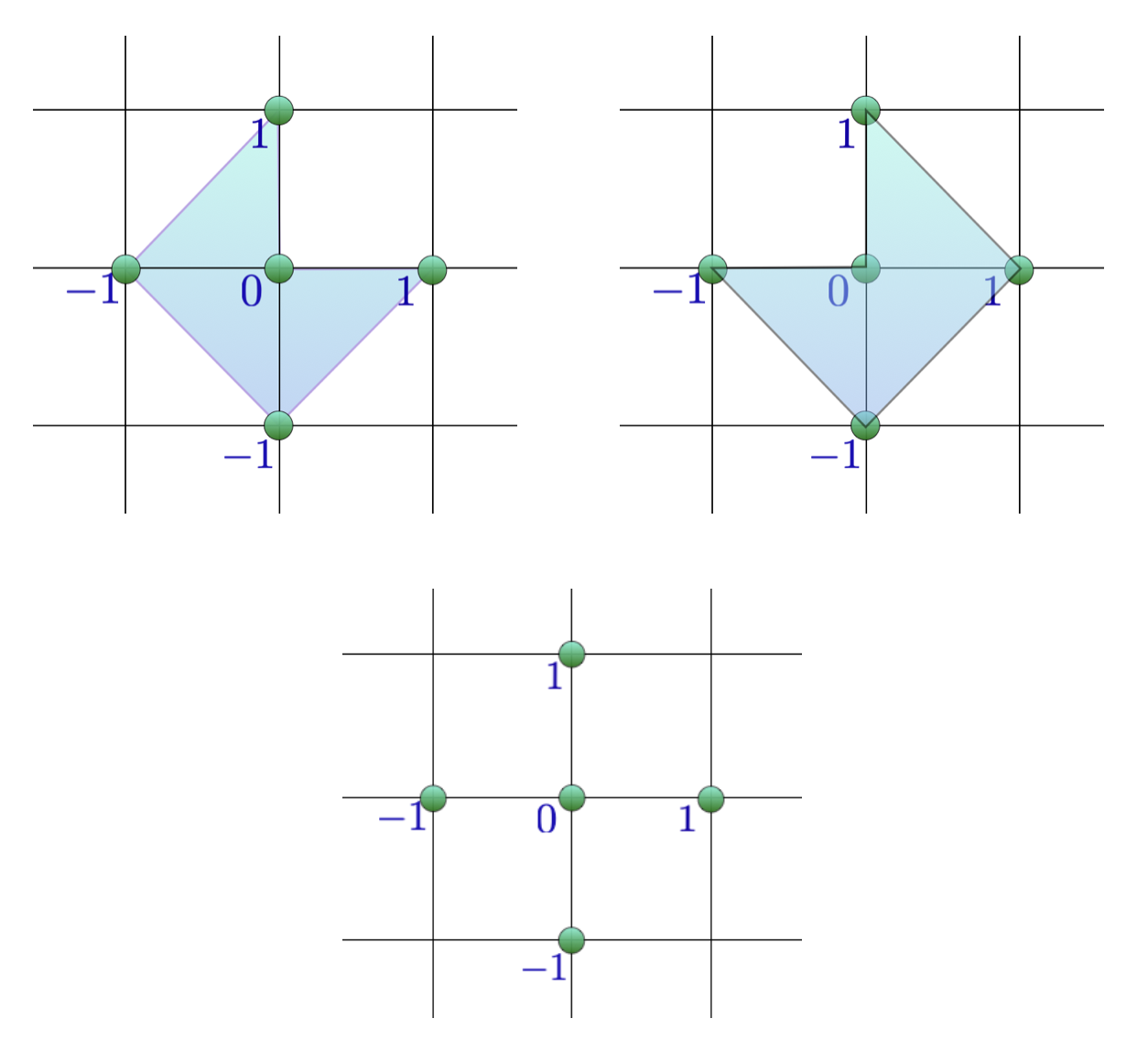}
		\end{center}
		\caption{Bottom: a~finite set of $5$ integer points. \ Top: two distinct, nonconvex polygons built on the same set of $5$ integer points. }
		\label{Nonconvex}
	\end{figure}

	Let us start with a~non-example, which will show the importance of convexity when passing from arbitrary sets of integer points in $\Z^d$ to the family of integer polytopes. Suppose we have two \emph{non-convex polytopes}
	$\PP, \QQ \subset \R^d$. It would be nice
	to say that if $\PP \cap \Z^d = \QQ \cap \Z^d$, then $\PP = \QQ$. Unfortunately, even in $\R^2$ we have simple counterexamples, as Figure~\ref{Nonconvex} shows.
	
	Fortunately, we have the very
	easy fact (Lemma~\ref{simple lemma, going from point set to polytope} below)
	that in the context of convex polytopes we do have such an implication, as follows.
	
	\begin{lemma}\label{simple lemma, going from point set to polytope}
		Let $\PP, \QQ \subset \R^d$ be two convex integer polyhedra.
		\begin{equation*}
			\text{If } \PP \cap \Z^d = \QQ \cap \Z^d, \text{ then } \PP = \QQ.
		\end{equation*}
	\end{lemma}
	\begin{proof}
		Each vertex of $\PP$ is an integer point, say $v\in \Z^d$, and is by assumption also a~point of $\QQ$. Therefore,
		the convex hull of the vertices (extreme points) of
		$\PP$ (which is $\PP$ itself, using the convexity of $\PP$) must be contained in $\QQ$, using the convexity of $\QQ$.
		So we have $\PP \subseteq \QQ$. By using an identical argument, we also have $\QQ \subseteq \PP$ and therefore $\PP = \QQ$.
	\end{proof}
	We will use Lemma~\ref{simple lemma, going from point set to polytope}
	repeatedly in the sequel, when moving from a~formal result about finite sets of integer points to results about polyhedra (and in particular polytopes).


	\section{Some examples}
	
	\begin{example}
		Let $\PP:= \conv\{(-1, 0), (1, 0), (0, 1)\}$. Its integer point transform is
		\[
		\sigma_\PP(\xi_1, \xi_2) = e^{2\pi i (-\xi_1)} + e^{2\pi i \xi_1} + e^{2\pi i \xi_2} \in \C.
		\]
		Evaluation at the special value
		$\xi^* :=~\left(\tfrac{\sqrt 2}{\pi}, \tfrac{\sqrt 3}{\pi}\right)$
		gives us a~unique `signature' for this polygon $\PP$. Namely, Theorem~\ref{first main theorem} tells us that
		$\PP$ is the only integer polygon associated to the special value
		$ \sigma_\PP\left(\tfrac{\sqrt 2}{\pi}, \tfrac{\sqrt 3}{\pi}\right)=
		e^{-i 2 \sqrt 2} +2e^{i 2 \sqrt 3}$.\hfill $\square$
	\end{example}

	\begin{figure}[htb]
		\begin{center}
			\includegraphics[totalheight=2.4in]{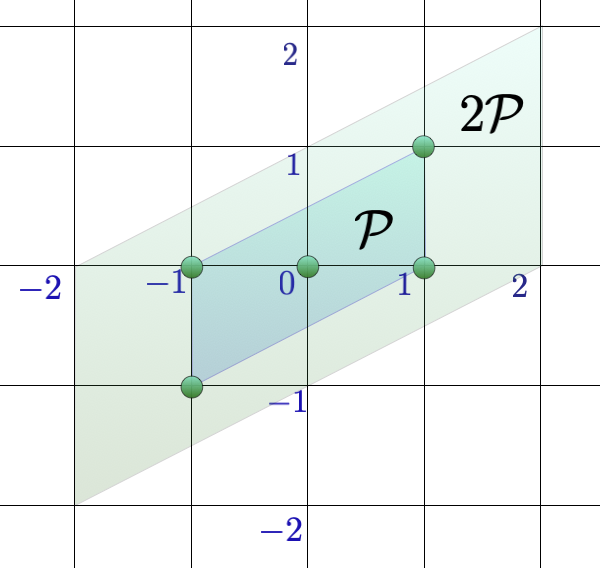}
		\end{center}
		\caption{
			The polygons $\PP$ and $2\PP$, of Example~\ref{symmetric 4-gon}
		}
		\label{Example1}
	\end{figure}

	\begin{example}\label{symmetric 4-gon}
		Let $\PP$ be the symmetric parallelogram whose vertices are given by
		\[
		\{ (1, 0), (1, 1), (-1, 0), (-1, -1)\},
		\]
		as in Figure~\ref{Example1}. Its integer point transform is
		\begin{align*}
			\sigma_\PP(\xi_1, \xi_2) &=1+ e^{2\pi i (\xi_1)} + e^{2\pi i (-\xi_1)} + e^{2\pi i (\xi_1 +\xi_2)} +e^{2\pi i (-\xi_1 -\xi_2)}  &\\
			&=1 + 2\cos(2\pi \xi_1) + 2 \cos(2\pi(\xi_1 + \xi_2)). & \square
		\end{align*}
	\end{example}
	
	\begin{example}\label{tetrahedron of index 2}
		Let $\PP$ be the integer tetrahedron in $\R^3$ whose vertices are given by
		\[
		(0, 0, 0), (1, 1, 0), (0, 1, 1), (1, 0, 1).
		\]
		The integer point transform of $\PP$ is
		\begin{equation*}
			\sigma_\PP(\xi)=1+e^{2\pi i (\xi_1 + \xi_2)}+e^{2\pi i (\xi_2 + \xi_3)}+e^{2\pi i (\xi_1 + \xi_3)}.
		\end{equation*}
		In Example~\ref{continuing with the tetrahedron of index 2}
		of Section~\ref{section: sublattices gen. by all the integer points in P}, we will use the same tetrahedron.
		
		The graph belies a~symmetry, as follows. If we replace any coordinate $\xi_k$ by $1-\xi_k$, the
		integer point transform in this example stays invariant. For example, when $k=1$:
		\begin{align*}
			\sigma_{\PP}(1-\xi_1, \xi_2, \xi_3, \cdots \xi_d)
			&= \sum_{n\in \Z^d \cap \PP}
			e^{2\pi i \big(n_1 (1-\xi_1) + n_2 \xi_2 + \cdots n_d \xi_d \big)} & \\
			&= \sum_{n\in \Z^d \cap \PP}
			e^{2\pi i \big(-n_1 \xi_1 + n_2 \xi_2 + \cdots n_d \xi_d \big)}. & \square
		\end{align*}
	\end{example}


	\section{The Fourier transform and integer point transform of a~general polyhedron, as complete invariants}\label{section:FT of integer stuff as a complete invariant}
	Here we extend Theorem~\ref{first main theorem} to show that the Fourier transform
	of a~rational polytope $\PP$ uniquely determines $\PP$ by evaluating the Fourier transform at a single algebraic point.
	Given a~$d$-dimensional compact set $S\subset \R^d$, we define the continuous
	Fourier transform of $S$ by
	\begin{equation}\label{Continuous FT}
		\mathcal F(S)(\xi):= \int_S e^{-2\pi i \langle \xi, x\rangle} dx,
	\end{equation}
	for each $\xi \in \R^d$. The present notation $\mathcal F$ is used to distinguish the Fourier transform~\eqref{Continuous FT} from the notation that we used for the Fourier transform over a~finite abelian group, namely $\mathcal F_G$.

	\begin{theorem}\label{thm:the FT as a complete invariant}
		We are given two rational polytopes $\PP, \QQ \subset \R^d$, and the algebraic
		point
		$\xi^*$ from~\eqref{alg. indep. vector}. Then we have:
		\begin{equation*}
			\mathcal F(1_\PP)(\xi^*) = \mathcal F(1_\QQ)(\xi^*)
			\iff
			\PP=\QQ.
		\end{equation*}
	\end{theorem}
	\begin{proof}
		We first recall Brion's theorem (\cite{Robins}, Theorem 8.3):
		\begin{equation*}
			\int_\PP e^{-2\pi i \langle u, \xi \rangle} \, du =
			\sum_{v \in V}
			\frac{e^{-2\pi i \langle v, \xi \rangle} }{(2\pi i)^d}
			\sum_{j=1}^{M(v)} \frac{\det \K_j(v) }{\prod_{k=1}^d \langle w_{j, k}(v), \xi \rangle},
		\end{equation*}
		for all $\xi \in \R^d$ such that none of the denominators vanish:
		$ \prod_{k=1}^d \langle w_{j, k}(v), \xi \rangle \not=0$.
		At each vertex $v \in \PP$, the vertex tangent cone $\K_v $
		is triangulated into simplicial cones, using the notation $\K_v = \K_1(v)\cup \dots \cup \K_{M(v)}(v)$.
		
		By assumption, all of the edge vectors $w_{j, k}(v)$ are integer vectors, the vertices $v$ are rational vectors (for all vertices of both $\PP$ and $\QQ$), and all the determinants $\det \K_j(v)$ are polynomial functions of the vertices $v$, and are therefore rational numbers.
		If we have $\mathcal F(1_\PP)(\xi^*) = \mathcal F(1_\QQ)(\xi^*) $, then
		\begin{equation*}
			\sum_{v \in V(\PP)}
			\frac{e^{-2\pi i \langle v, \xi^* \rangle} }{(2\pi i)^d}
			\sum_{j=1}^{M(v)} \frac{\det \K_j(v) }{\prod_{k=1}^d \langle w_{j, k}(v), \xi^* \rangle}
			=
			\sum_{v \in V(\QQ)}
			\frac{e^{-2\pi i \langle v, \xi^* \rangle} }{(2\pi i)^d}
			\sum_{j=1}^{M(v)} \frac{\det \K_j(v) }{\prod_{k=1}^d \langle w_{j, k}(v), \xi^* \rangle}
		\end{equation*}
		holds, and this may be rewritten as
		\begin{equation}\label{the exponential-polynomiall setup}
			0=
			\sum_{v \in V(\PP)} \mathrm{c_\PP}(v, \xi^*)
			e^{-2\pi i \langle v, \xi^* \rangle}
			-
			\sum_{v \in V(\QQ)} \mathrm{c_\QQ}(v, \xi^*)
			e^{-2\pi i \langle v, \xi^* \rangle},
		\end{equation}
		where the coefficients $\mathrm{c_\PP}(v, \xi^*)$ and $\mathrm{c_\QQ}(v, \xi^*)$ are polynomial functions of the coordinates of the rational vertices $v$, with algebraic coefficients, due to the appearance of $\xi^*$.
		Supposing that $\PP \not=\QQ$, we have arrived at a~nontrivial linear combination of
		exponentials of the form $e^{\alpha_1}, \dots, e^{\alpha_n}$, with algebraic coefficients, and with
		$\alpha_1, \dots, \alpha_n$ distinct algebraic numbers, by our choice of $\xi^*$.
		But this contradicts the Lindemann–Weierstrass theorem
		(\cite{Baker}, Theorem 1.4).
	\end{proof}
	

	\section{Integer sublattices generated by the integer points in a~polytope, and absolute maxima of integer point transforms}\label{section: sublattices gen. by all the integer points in P}
	
	Let $\PP\subset \R^d$ be an integer polytope. We define
	\[
	\Lat_\PP\subset \Z^d
	\]
	to be the integer sublattice of $\Z^d$ that is generated by the integer span of all
	the points in $| \Z^d \cap \PP|$. We'll call $ \Lat_\PP$ the {\bf spanning lattice} of the set $\PP\cap \Z^d$.  It is a~trivial fact\textemdash perhaps as old as the hills themselves\textemdash that for $2$-dimensional integer polygons, we always have $\Lat_\PP = \Z^2$. The proof is easy: we can always triangulate any integer polygon into unimodular triangles, and the vertices of any such unimodular triangle already generate $\Z^2$.
	
	Whenever $\Lat_\PP= \Z^d$, $\PP$ is called a~\emph{spanning polytope}. In dimensions $d\geq 3$, there are integer polytopes that are not spanning polytopes, because a~unimodular triangulation is not always available for an arbitrary integer polytope.
	\begin{example}\label{Reeve tetrahedron}
		Given any positive integer $h$, the \emph{Reeve} tetrahedron $T_h \subset \R^3$ has vertices
		$\left(0, 0, 0\right), \left(1, 0, 0\right), \left(0, 1, 0\right), \left(1, 1, h\right)$.
		It is an exercise that $T_h$ is not a~spanning polytope, for any $h \geq 2$.
		\hfill $\square$
	\end{example}
	In the recent work~\cite{Hofs}, Hofscheier, Katth\"an, and Nill develop an Ehrhart-type theory for spanning lattice polytopes.
	Here we may ask another natural question about the integer point transform:
	\begin{question}\label{absolute maxima}
		Which values of $\xi \in \R^d$ give us the absolute maxima of $\left | \sigma_\PP(\xi) \right |$?
	\end{question}
	It turns out that the dual lattice $\Lat_\PP^*$ comes in naturally here, and it satisfactorily answers
	Question~\ref{absolute maxima}.

	\begin{theorem}\label{thm:absolute maxima}
		For an integer polytope $\PP\subset \R^d$ that contains the origin, we have:
		\[
		\left | \sigma_\PP(\xi) \right | = \left | \Z^d \cap \PP \right | \iff \xi \in \Lat_\PP^*, 
		\]
		the dual of the spanning lattice. 
	\end{theorem}
	\begin{proof}
		The triangle inequality for complex numbers tells us that
		\begin{equation}\label{lemma:triangle inequality}
			\left | \sigma_{\PP}(\xi) \right |
			\leq
			\sum_{n\in \Z^d \cap \PP} \left | e^{2\pi i \langle n, \xi \rangle} \right |
			= \sum_{n\in \Z^d \cap \PP}1 = \left | \Z^d \cap \PP \right |,
		\end{equation}
		for all $\xi \in \R^d$. Let us fix $\xi \in \R^d$ such that $\left | \sigma_\PP(\xi) \right | = \left | \Z^d \cap \PP \right |$.
		By~\eqref{lemma:triangle inequality}, the latter equality means
		\begin{equation}\label{modsumeq}
			\left | \sum_{n\in \Z^d \cap \PP} e^{2\pi i \langle n, \xi \rangle} \right |
			=\sum_{n\in \Z^d \cap \PP} \left | e^{2\pi i \langle n, \xi \rangle} \right |,
		\end{equation}
		
		which in turn occurs exactly when all of the complex numbers $e^{2\pi i \langle n, \xi \rangle}$ point in the same direction:
		\begin{align}
			\eqref{modsumeq} & \iff e^{2\pi i \langle n, \xi \rangle}=e^{2\pi i \langle x_0, \xi \rangle} \ \text{ for some } x_0 \in \PP,
			\text{ and for all } n \in \Z^d \cap \PP    \\
			& \iff \langle n-x_0, \xi \rangle \in \Z \text{ for all } n \in \Z^d \cap \PP. \label{almost at the dual lattice}
		\end{align}
		
		We recall that by definition $\Lat_\PP$ is generated by the integer span of all $n \in \Z^d \cap \PP$. So
		the condition~\eqref{almost at the dual lattice} holds if and only if $\langle n, \xi \rangle \in \Z $ for all $n \in \Lat_\PP$, by the linearity of the inner product. By definition, this means that $\xi \in \Lat_\PP^*$, the dual lattice.
	\end{proof}
	
	We notice that if replace the polytope $\PP$ by any \emph{finite subset of integer points}
	$S\subset\Z^d$, then the proof of
	Theorem~\ref{thm:absolute maxima} remains valid, with the following definition. We define
	$\Lat_{S}\subset \Z^d$
	to be the lattice generated by the integer span of all $n \in S$. We obtain the following immediate consequence from the proof of Theorem~\ref{thm:absolute maxima}.

	\begin{corollary}\label{cor: absolute maxima equivalence}
		Given any finite set of integer points $S\subset \Z^d$, we have
		\begin{enumerate}
			\item $\left | \sigma_S(\xi) \right | \leq |S|$ for all $\xi \in \R^d$.
			\item $\left | \sigma_S(\xi) \right | = |S| \iff \xi \in \Lat_S^*$.
			\label{cor:part 2 of our abs. maxima equivalence}
		\end{enumerate}
	\end{corollary}
	
	Using the definition of a~spanning polytope, together with Theorem~\ref{thm:absolute maxima}, we immediately obtain the following consequence as well.
	\begin{corollary}
		For an integer polytope $\PP\subset \Z^d$, the following are equivalent:
		\begin{enumerate}
			\item $\PP$ is a~spanning polytope.
			\item
			$
			\left | \sigma_\PP(\xi) \right | = \left | \Z^d \cap \PP \right | \iff \xi \in \Z^d.
			$
		\end{enumerate}
	\end{corollary}

	Given the nice structure of the aboslute maxima given by Theorem~\ref{thm:absolute maxima}, it is also natural to ask:
	\begin{question}\label{abs. maxima modulo the unit cube}
		How many inequivalent absolute maxima are there, modulo $\Z^d$?
	\end{question}
	We have the sublattice containments
	$\Lat_{\PP} \subseteq \Z^d \subseteq \Lat_\PP^*$.
	Question~\ref{abs. maxima modulo the unit cube} asks for the value of $[\Lat_\PP^* : \Z^d ]$,
	which has a~simple answer:
	\begin{equation}\label{number of inequivalent maxima}
		[\Lat_\PP^* : \Z^d ] = \frac{ \det \Z^d}{\det \Lat_\PP^*} = \frac{ 1}{\left(1/\det \Lat_\PP\right)}
		= \det \Lat_\PP.
	\end{equation}

	\begin{example}\label{continuing with the tetrahedron of index 2}
		Let us use Theorem~\ref{thm:absolute maxima} to find all of the absolute maxima of the
		integer point transform for the tetrahedron in
		Example~\ref{tetrahedron of index 2}. We recall that
		\[
		\PP:= \conv\{(0, 0, 0), (1, 1, 0), (0, 1, 1), (1, 0, 1)\},
		\]
		had
		$
		\sigma_\PP(\xi)=1+e^{2\pi i (\xi_1 + \xi_2)}+e^{2\pi i (\xi_2 + \xi_3)}+e^{2\pi i (\xi_1 + \xi_3)}.
		$
		Here the integer span of the $4$ integer points that comprise $\PP\cap \Z^3$ gives us the sublattice
		$\Lat_\PP:= M(\Z^3)$, with
		\[
		M:=
		\begin{pmatrix}
			1 & 1 & 0\\
			1 & 0 & 1 \\
			
			0 & 1 & 1
		\end{pmatrix}
		.
		\]
		This is a~sublattice of index $2$ in $\Z^3$, because $\det c = \det M=2$. Here the dual lattice $\Lat^*$ has a~generator matrix:
		\[
		M^{-T}:=
		\begin{pmatrix}
			\ \ \ \frac{1}{2} & \ \ \ \frac{1}{2} & -\frac{1}{2} \\
			\ \ \ & \ \ \ & \\
			\ \ \ \frac{1}{2} & \ -\frac{1}{2} & \ \ \frac{1}{2} \\
			\ \ \ & \ \ \ & \\
			-\frac{1}{2} & \ \ \frac{1}{2}& \ \ \frac{1}{2}
		\end{pmatrix}
		.
		\]
		Let us check Corollary~\ref{cor: absolute maxima equivalence},
		part~\ref{cor:part 2 of our abs. maxima equivalence} concerning the locations of absolute maxima. First, one of the absolute maxima of $\sigma_\PP$ is $|\sigma_\PP(0)| = 4$.
		Now, we have
		\[
		v:=
		\begin{pmatrix}
			\frac{1}{2} \\
			\\
			\frac{1}{2} \\
			\\
			\frac{1}{2}
		\end{pmatrix}
		\in \Lat_\PP^*,
		\]
		because $v$ equals the sum of the three columns of $M^{-T}$.
		Let us check if $v$ gives us another absolute maxima, as predicted by
		Corollary~\ref{cor: absolute maxima equivalence}, part~\ref{cor:part 2 of our abs. maxima equivalence}:
		\[
		|\sigma_\PP(v)| = \left |
		1+e^{2\pi i \left(\tfrac{1}{2} + \tfrac{1}{2}\right) }+e^{2\pi i \left(\tfrac{1}{2} - \tfrac{1}{2}\right)}
		+e^{2\pi i \left(\tfrac{1}{2} - \tfrac{1}{2}\right)}
		\right|
		=4,
		\]
		so indeed it does. Due to the relation~\eqref{number of inequivalent maxima}, we know that the total number of inequivalent absolute maxima in $[0, 1)^d$ is equal to
		$\det \Lat_\PP = 2$, so we found both of them.
		\hfill $\square$
	\end{example}


	\section{Integer point transforms on a~finite abelian group}\label{machinery of characters on finite abelian groups}
	
	The finite sum of exponentials that defines $\sigma_\PP(\xi)$ lends the feeling that we should be studying a~connection to the finite Fourier transform of some finite abelian group. To make this feeling rigorous, we develop this connection here. Although there are many possible choices for our finite abelian group, we first choose a~box
	\begin{equation}\label{the box}
		B:=
		\left[-\frac{k_1}{2}, \frac{k_1}{2}\right] \times \cdots \times \left[-\frac{k_d}{2}, \frac{k_d}{2}\right],
	\end{equation}
	and then consider the set of integers points in it, namely $\Z^d \cap B$.
	Clearly, any integer polytope $\PP$ is contained in the interior of the box $B$, for some appropriately chosen integers $k_1, \dots, k_d$. To avoid ambiguities when we take the quotient, we assume that each vertex of $\PP$ is
	contained in the slightly smaller box
	$[-\frac{k_1}{2} +\frac{1}{2}, \frac{k_1}{2}-\frac{1}{2}] \times \cdots \times
	[-\frac{k_d}{2} +\frac{1}{2}, \frac{k_d}{2}-\frac{1}{2}]$.

	Now we consider the finite abelian group $G:= \Z / k_1 \Z \times \cdots \Z/ k_d \Z$, which we identify with the integer points in the torus $\R / k_1 \Z \times \cdots \times \R / k_d \Z$, and which can also be thought of
	as the box $B$ after identifying its opposite facets.
	Although the choice of integers $k_1, \dots k_d$ is not canonical,
	such an embedding of the integer points of $\PP$ into a~finite abelian group will prove to be worthwhile. In other words, we have, by definition:
	\[
	\PP \cap \Z^d := \text{ the domain of } 1_P, \text{ as a~function on } G.
	\]
	We will freely use the usual fact that the Pontryagin dual $\hat G$ (simply the group of characters of $G$) is in this case isomorphic to $G$.
	Even though the full geometry of $\PP$ may not be immediately
	apparent in this discrete setting, we will be able to shed some additional light on the integer points
	$
	\PP \cap \Z^d \subset B.
	$
	It is now natural to consider the indicator function $1_\PP$ as a~function on $G$,
	and therefore expand it into
	its finite Fourier series. Precisely, each element $\xi \in G$ gives us a~character $\chi_\xi: G\rightarrow S^1$
	defined by
	\[
	\chi_\xi(n):=e^{2\pi i \left(\frac{\xi_1n_1}{k_1}+\cdots+\frac{\xi_dn_d}{k_d}\right)}.
	\]
	The first theorem of finite Fourier analysis gives us:
	\begin{equation*}
		1_\PP(n) = \sum_{\xi \in G} c_\xi \chi_\xi(n),
	\end{equation*}
	where the (finite) Fourier coefficients have the form
	\begin{equation*}
		c_\xi = \sum_{g \in G} 1_\PP(g) \overline{ \chi_\xi (g) },
	\end{equation*}
	and where we have used the isomorphism $G \cong \hat G$.
	So by definition we have the \emph{finite Fourier transform}
	\[
	\mathcal F_G(1_P)(\xi):= c_\xi.
	\]
	Let us ``massage'' $c_\xi$ a~bit. For each $\xi \in G$, we have:
	\begin{align}
		\mathcal F_G(1_P)(\xi)= c_\xi &= \sum_{g \in G} 1_\PP(g) \overline{ \chi_\xi(g) } \\
		&= \sum_{g \in G} 1_\PP(g) e^{-2\pi i \left(\frac{g_1 \xi_1}{k_1}+\cdots+\frac{g_d \xi_d}{k_d}\right)} \\
		&= \sum_{m \in \PP \cap \Z^d} e^{-2\pi i \left(\frac{m_1 \xi_1}{k_1}+\cdots+\frac{m_d\xi_d}{k_d}\right)}\\
		\label{identifying the finite FT with the integer point transform}
		&:= \sigma_\PP\left(-\frac{\xi_1}{k_1}, \dots, -\frac{\xi_d}{k_d}\right),
	\end{align}
	by definition of the integer point transform.
	We will also use the latter identification in Section~\ref{sec:central symmetry}
	below, by giving an equivalent condition for central symmetry in terms of finite Fourier transforms\textemdash or equivalently the integer point transform.
	
	We now notice that we have never required $\PP$ to be a~polytope in the theory above,
	but merely that we have:
	\[
	\text{\emph{any finite subset of the integer lattice}}.
	\]
	One reason for initially using the integer points that belong to a polytope, as opposed to just any finite set of integer points, is that the applications that use polytopes are perhaps the most naturally occurring.
	
	Next, suppose that we want an analogue of Theorem \ref{first main theorem}, but we want to evaluate the integer point transform at integer lattice points, rather than evaluating it at the algebraic point $\xi^* := \tfrac{1}{\pi}\left(\sqrt 2, \dots, \sqrt{p_d}\right)$. Then we need to evaluate at more points, and the next result, namely
	Theorem~\ref{Finite Fourier transform evaluation}, part~\ref{part b of Finite Fourier transform evaluation}, gives a~sufficient condition to choose such points.

	\begin{theorem}\label{Finite Fourier transform evaluation}
		Suppose that $S\subset \Z^d$ is a~finite subset of the half-open box
		\[
		B:=
		\left[-\frac{k_1}{2}, \frac{k_1}{2}\right) \times \cdots \times \left[-\frac{k_d}{2}, \frac{k_d}{2}\right).
		\]
		With the notation above, the following hold:
		%
		%
		\begin{enumerate}
			\item \label{part a of Finite Fourier transform evaluation}
			\begin{align*}
				\mathcal F_G(1_S)(\xi)=\sigma_S\left(-\frac{\xi_1}{k_1}, \dots, -\frac{\xi_d}{k_d}\right).
			\end{align*}
			\item \label{part b of Finite Fourier transform evaluation}
			$S$ is uniquely determined by the finite
			set of special values
			\[
			\left\{
			\sigma_S\left(\frac{\xi_1}{k_1}, \dots, \frac{\xi_d}{k_d}\right) \mid \
			\xi \in B\cap \Z^d
			\right\}.
			\]
		\end{enumerate}
	\end{theorem}
	\begin{proof}
		Part~\ref{part a of Finite Fourier transform evaluation} follows from the definition of $\sigma_S$ and our discussion of finite Fourier transforms above. For part~\ref{part b of Finite Fourier transform evaluation}, we will use
		the uniqueness of the inverse Fourier transform over the finite abelian group $G$ defined above. In particular, suppose that we have two sets $S_1, S_2 \subset B$, or equivalently $S_1, S_2 \subset G$. Then we have $\mathcal F_G(1_{S_1})(\xi) =\mathcal F_G(1_{S_2}(\xi)$ for all $\xi \in G$. By Fourier inversion, we have
		\begin{align*}
			1_{S_1}(\xi)
			&=\mathcal F_G \left(\mathcal F_G(1_{S_1})\right)(-\xi)
			= \mathcal F_G \left(\mathcal F_G(1_{S_2})\right)(-\xi) \\
			&=1_{S_2}(\xi),
		\end{align*}
		for all $\xi \in G$. Therefore $1_{S_1}=1_{S_2}$, and so $S_1 = S_2$.
	\end{proof}


	\section{Centrally symmetric sets of integer points and centrally symmetric polytopes}\label{sec:central symmetry}
	
	In this brief section we give an equivalence for central symmetry in terms of the integer point transform, for any finite set
	$A\subset \Z^d$ of integer points. As a~consequence, we get an equivalence for the central symmetry of any integer polytope (of arbitrary codimension) in terms of special evaluations of the integer point transform. We will use here the machinery of Section~\ref{machinery of characters on finite abelian groups}.
	
	In this section we will use the box
	$B:=
	\left[-\frac{k}{2}, \frac{k}{2}\right)^d$, for any positive integer $k$. We will suppose that
	$A\subset \Z^d \cap B$ is any centrally symmetric set of integer points that is contained in $B$.
	
	It is immediate that for such a~set $A$, its integer point transform is real-valued for any $\xi \in \R^d$:
	\begin{align}
		\sigma_{A}(\xi):= \sum_{n\in A} e^{2\pi i \langle n, \xi \rangle}
		&:= 1 + \frac{1}{2} \sum_{\genfrac{}{}{0pt}{}{m\in A}{m \not=0}} e^{2\pi i \langle m, \xi \rangle}
		+ \frac{1}{2} \sum_{\genfrac{}{}{0pt}{}{m\in -A}{m \not=0}} e^{2\pi i \langle m, \xi \rangle}\\
		&:= 1 + \frac{1}{2} \sum_{\genfrac{}{}{0pt}{}{m\in A}{m \not=0}} e^{2\pi i \langle m, \xi \rangle}
		+ \frac{1}{2} \sum_{\genfrac{}{}{0pt}{}{m\in A}{m \not=0}} e^{2\pi i \langle -m, \xi \rangle}\\
		&= 1+ \sum_{\genfrac{}{}{0pt}{}{m\in A}{m \not=0}} \cos\left(2\pi \langle m, \xi \rangle\right) \in \R.
		\label{one direction of cs polytope giving real transform}
	\end{align}
	In particular, for centrally symmetric polytopes, their minima and maxima may be studied without
	having to take norms.

	\begin{theorem}\label{thm:central symmetry equivalence}
		Let $A\subset \Z^d$ be a~finite collection of integer points.
		The following are equivalent:
		\begin{enumerate}
			\item $A$ is centrally symmetric.
			\label{part 1 of central symmetry}
			\item $\sigma_A(\frac{1}{k}\xi) \in \R$ for all $\xi \in A$.
			\label{part 2 of central symmetry}
		\end{enumerate}
	\end{theorem}
	\begin{proof}
		For the easy direction that~\eqref{part 1 of central symmetry} $\implies$~\eqref{part 2 of central symmetry}, we have already seen the proof
		in~\eqref{one direction of cs polytope giving real transform}. Now suppose that
		$\sigma_A(\frac{1}{k}\xi) \in \R$ for all $\xi \in A$, and we must show that $A = -A$.
		We therefore have
		$\sigma_A(\frac{1}{k}\xi) = \overline{\sigma_A(\frac{1}{k}\xi)}$, and using
		equation~\eqref{identifying the finite FT with the integer point transform} (with all $k_j= k$)
		we may rewrite this condition in terms of the finite Fourier transform as
		\begin{align}
			\mathcal F_G(1_A)\left(\xi\right) \label{penultimate identity in the proof of central symmetry}
			&= \overline{\mathcal F_G(1_A)}\left(\xi\right) \\
			&= \sum_{m \in A} e^{2\pi i \left(\frac{m_1 \xi_1}{k}+\cdots+\frac{m_d\xi_d}{k}\right)}\\
			&= \sum_{m \in (-A)} e^{-2\pi i \left(\frac{m_1 \xi_1}{k}+\cdots+\frac{m_d\xi_d}{k}\right)}\\
			&=\mathcal F_G(1_{\{-A\}})\left(\xi\right),
		\end{align}
		for all $\xi \in A$. Finally, we now take the inverse Fourier transform of both sides, to conclude that
		$1_A\left(\xi\right) =1_{\{-A\}}\left(\xi\right)$ for all $\xi \in A$.
		Therefore $A = -A$.
	\end{proof}

	
	\section{Integer point transforms and Fourier transforms of integer cones are also complete invariants}
	
	Using exactly the same proof ideas of Theorem \ref{first main theorem} and 
	Theorem \ref{thm:the FT as a complete invariant}, we also obtain the following corollaries.
	
	\begin{corollary}\label{thm:the FT of a cone as a complete invariant}
		Given any two integer cones $\K_1, \K_2 \subset \R^d$, we have:
		\begin{equation*}
			\mathcal F(1_{\K_1})(\xi^*) = \mathcal F(1_{\K_2})(\xi^*)
			\iff
			\K_1=\K_2,
		\end{equation*}
		where $F(1_{\K})$ is the Fourier-Laplace transform of the cone $\K$.
	\end{corollary}
	\begin{proof} A standard and known computation \cite[Corollary 8.1]{Robins} 
		gives us the Fourier transform of a cone $\K$:
		\begin{equation}\label{FT of ANY cone}
			F(1_{\K})(\xi) = 
			\int_{\K} e^{-2\pi i \langle u, \xi \rangle} \, du =
			\sum_{j=1}^{M}
			\frac{e^{-2\pi i \langle v, \xi \rangle} }{(-2\pi i)^d}
			\frac{\det \K_j}{\prod_{k=1}^d \langle w_{j, k}, \xi \rangle},
		\end{equation}
		for all $\xi \in \R^d$ such that none of the denominators vanish. Here we've 
		triangulated the cone $\K$ into simplicial cones $\K_1, \dots, \K_M$. 
		
		We note that, initially,
		formula \eqref{FT of ANY cone} holds for a complex vector $\xi$ that allows the defining integral over $\K$ to converge.  However, by meromorphic continuation, we may later plug in any real $\xi \in \R^d$, as long as the denominators in  \eqref{FT of ANY cone} do not vanish. 
		The rest of the proof is identical to the proof of Theorem  \ref{thm:the FT as a complete invariant}.
	\end{proof}
	
	\begin{corollary}\label{cor:integer cones are also uniquely determined}
		Given any two integer cones $\K_1, \K_2 \subset \R^d$,  we have
		\begin{equation*}
			\sigma_{\K_1}(\xi^*) =\sigma_{\K_2}(\xi^*) \iff \K_1 = \K_2.
		\end{equation*}
	\end{corollary}
	\begin{proof}
		We recall an elementary lemma (for example~\cite[Theorem 10.2]{Robins}), which tells us that 
		the integer point transform $\sigma_{\K}$ of  a~simplicial integer cone $\K$ has the following finite form, as a rational function:
		\begin{equation}\label{rational function setup for cones}
			\sigma_{\K}(\xi)
			=
			\frac{ \sigma_{ \Pi }(\xi) }
			{ \prod_{k=1}^d \left(1 - e^{ \langle w_k, \xi \rangle}\right) },
		\end{equation}
		where $\Pi := \{ \lambda_1 w_1 + \cdots + \lambda_d w_d  \,
		\mid \text{ all } 0 \leq \lambda_j < 1 \}$, and where $w_1, \dots, w_d\in \Z^d$ are the integer edge vectors of $\K$.  Here,  $\sigma_{ \Pi }(\xi)$ is the integer point transform of a finite set of integer points, and hence a (Laurent) polynomial. For any (possibly non-simplicial) integer cone $\K$, it is also a fact that its integer point transform is a finite linear combination over $\Z$, of rational functions identical to \eqref{rational function setup for cones}.
		
		If $\sigma_{\K_1}(\xi^*) =\sigma_{\K_2}(\xi^*)$,  then using 
		\eqref{rational function setup for cones}
		we arrive at an identical equation to 
		\eqref{the exponential-polynomiall setup}, and the rest of the proof is identical to the proof of
		Theorem \ref{thm:the FT as a complete invariant}.
	\end{proof}

	\section{Further remarks and questions}\label{Further remarks}
	For further rumination, we mention a~few threads that appeared naturally in this line of research which remain open.

	\begin{enumerate}
		\item
		\label{first remark}
		Perhaps the most fascinating question now is how to find the unique polytope (or set of integer points) that is guaranteed by the uniqueness property of
		Theorem~\ref{first main theorem}.
		
		For example, suppose we seek to discover the
		$1$-dimensional polytope $\PP:= [0, N]$, and suppose its integer point transform
		$\sigma_\PP(\xi^*) = C$ is given to us, with $\xi^*$ as in~\eqref{alg. indep. vector}. Then we have
		\[
		C = \sum_{k =0}^N e^{2\pi i k \xi^* } = \frac{ e^{2\pi i (N+1) \xi^* }-1 }{e^{2\pi i \xi^* }-1}.
		\]
		Here, it is easy to solve this equation with respect to the vertex $N$ of $\PP$, by taking complex logs, with some care being taken for picking the principal branch. However, even for an arbitrary finite set $\PP$ of integer points in $\Z$, or for example an integer triangle in $\R^2$, the problem already appears to become formidable.
		
		\item
		\label{second remark}
		Similarly, it would be important to reconstruct a~rational polytope from the uniqueness guaranteed by the 1-point evaluation of its continuous Fourier transform, as in Theorem~\ref{thm:the FT as a complete invariant}. This direction for future research, as well as the previous problem, involves transcendental equations. It would be very interesting to solve them over the integers, namely the coordinates of the vertices of $\PP$.
		
		A similar open problem is to solve for the unique integer cone that Corollary 
		\ref{thm:the FT of a cone as a complete invariant} and Corollary
		\ref{cor:integer cones are also uniquely determined} guarantee, respectively.

		\item \label{third remark}
		Is it possible to strengthen Theorem~\ref{first main theorem}, part~\ref{second part of main theorem}, so that we can eliminate the dilation factor $k$, as follows?
		\begin{conj}
			If $\PP, \QQ \subset \R^d$ are any rational polytopes, then we have:
			\[
			\sigma_{\PP}(\xi^*) =\sigma_{\QQ}(\xi^*) \implies \PP = \QQ,
			\]
			with $\xi^*$ as in~\eqref{alg. indep. vector}. \hfill $\square$
		\end{conj}

		\item
		\label{fourth remark}
		Historically, the integer point transform $\sigma_\PP$ for a~rational polytope $\PP$ appeared in the work of Brion~\cite{Brion} in 1988, who proved the important result that we may write
		$\sigma_\PP$ as a~finite linear combination of exponential-rational functions of the
		vertex tangent cones $\K_v$ of $\PP$:
		\begin{equation}\label{Brion's result}
			\sigma_\PP(\xi) = \sum_{\mathrm{ vertex } \, v \, { \rm of } \, \PP}
			\sigma_{\K_v}(\xi),
		\end{equation}
		valid for almost all $\xi \in \R^d$ (see~\cite{BeckRobins},~\cite{Robins} for more details). We have slightly abused notation in~\eqref{Brion's result} by writing $\sigma_{\K_v}(\xi)$ to mean the meromorphic continuation of these integer point transforms, which are initially defined by
		$\sigma_{\K_v}(\xi) := \sum_{n \in \K_v \cap \Z^d} e^{2\pi i \langle \xi, n \rangle}$.
		By an elementary lemma (for example~\cite[Theorem 10.2]{Robins}), each such $\sigma_{\K_v}$ for a~simplicial cone $\K$ also has the following finite form, as a rational-exponential function:
		\[
		\sigma_{\K_v}(\xi)
		=
		\frac{ \sigma_{ \Pi + v}(\xi) }
		{ \prod_{k=1}^d \left(1 - e^{ \langle w_k, \xi \rangle}\right) },
		\]
		where $\Pi +v := \{ \lambda_1 w_1 + \cdots + \lambda_d w_d +v \,
		\mid \text{ all } 0 \leq \lambda_j < 1 \}$, and where $w_1, \dots, w_d$ are the edge vectors of the vertex tangent cone $\K_v$.

		\item
		\label{fifth remark}
		Regarding possible extensions, it is tempting extend the integer point transform to arbitrary lattices, as follows.
		Suppose we are given any full-rank lattice $\Lat := M(\Z^d)$, and we define
		the integer point transform of a~given rational polytope $\PP \subset \R^d$, relative to $\Lat$,
		by
		\begin{equation*}
			\sigma_{\PP, \Lat}(\xi):= \sum_{n\in \Lat \cap \PP} e^{2\pi i \langle n, \xi \rangle},
		\end{equation*}
		for all $\xi \in \R^d$.
		Offhand, it may seem like we have a~new extension, but in fact we may easily rewrite it
		as follows:
		\begin{align*}
			\sigma_{\PP, \Lat}(\xi)&:=
			\sum_{n\in M(\Z^d) \cap \PP} e^{2\pi i \langle n, \xi \rangle}
			= \sum_{k\in \Z^d \cap \PP} e^{2\pi i \langle Mk, \xi \rangle}
			= \sum_{k\in \Z^d \cap \PP} e^{2\pi i \langle k, M^T \xi \rangle}
			:= \sigma_{\PP}(M^T \xi).
		\end{align*}
		Since $M$ is invertible ($\Lat$ has full-rank), and $\xi$ varies over all of $\R^d$, there is nothing really new in this particular extension, so we may use the usual integer point transform to sum over any lattice.
	\end{enumerate}


	{\small\bibliography{commat}}

	\EditInfo{April 19, 2023}{May 21, 2023}{Camilla Hollanti and Lenny Fukshansky}
	
\end{document}